\documentclass{amsart}
\usepackage{amssymb}
\usepackage{latexsym}
\usepackage{amsmath}
\usepackage{euscript}
\usepackage{graphics}

      \def\dC{{\mathbb C}}

   \def\dN{{\mathbb N}}   
      \def\dR{{\mathbb R}}

   \def\dZ{{\mathbb Z}}

\def\cJ{{\mathcal J}}

   \def\cW{{\mathcal W}}

\def\t#1{{{\tilde #1} }}

\def\bm\chi{\mbox{\boldmath$\chi$}}

\let\xker=\ker \def\ker{{\xker\,}}

\def\f{\varphi}

\def\nk{{\mathbf N}_\kappa}

\def\Im{\operatorname{Im}}

\def\dist{\operatorname{dist}}

\newcommand{\be}{\begin{equation}}

\newcommand{\ee}{\end{equation}}

\newcommand{\ba}{\begin{eqnarray}}
\newcommand{\ea}{\end{eqnarray}}

\newcommand{\baa}{\begin{eqnarray*}}
\newcommand{\eaa}{\end{eqnarray*}}
\newcommand{\bb}{\begin{thebibliography}}
\newcommand{\eb}{\end{thebibliography}}

\newcommand{\bi}[1]{\bibitem{#1}}
\newcommand{\lab}[1]{\label{#1}}
\newcommand{\re}[1]{(\ref{#1})}

\newtheorem{theorem}{Theorem}[section]
\newtheorem{proposition}[theorem]{Proposition}
\newtheorem{corollary}[theorem]{Corollary}
\newtheorem{lemma}[theorem]{Lemma}
\newtheorem{definition}[theorem]{Definition}
\theoremstyle{definition}

\newtheorem{remark}[theorem]{Remark}

\numberwithin{equation}{section}

\begin{document}

\title[An operator approach to multipoint Pad\'e approximations]
{An operator approach to multipoint Pad\'e approximations}

\author{Maxim S. Derevyagin}
\author{Alexei S. Zhedanov}

\address{Department of Nonlinear Analysis \\
Institute of Applied Mathematics and Mechanics\\
R.Luxemburg str. 74 \\
83114 Donetsk, Ukraine \\} \email{derevyagin.m@gmail.com}

\address{Institute for Physics and Engineering\\
R.Luxemburg str. 72 \\
83114 Donetsk, Ukraine \\}

\email{zhedanov@yahoo.com}

\date{\today}

\begin{abstract}
First, an abstract scheme of constructing biorthogonal rational
systems related to some interpolation problems is proposed. We
also present a modification of the famous step-by-step process of
solving the Nevanlinna-Pick problems for Nevanlinna functions. The
process in question gives rise to three-term recurrence relations
with coefficients depending on the spectral parameter. These
relations can be rewritten in the matrix form by means of two
Jacobi matrices. As a result, a convergence theorem for multipoint
Pad\'e approximants to Nevanlinna functions is proved.
\end{abstract}

\keywords{Nevanlinna-Pick problem, biorthogonal rational
functions, $R_{II}$-fraction, multipoint Pad\'e approximant,
Jacobi matrix, generalized eigenvalue problem, Markov function}

\maketitle

%%%%%%%%%%%%%%%%%%%%%%%%%%%%%%%%%%%%%%%%%%%%%%%%%%%%%%%%%%%%%%%%%%%%%%%%%%%%%%%%%%%%%%%%%%%%%%%%%%%%%%%%%%%%%%%%%%%%%%%%%%%%%%%%%%%%%%%%5

\section{Introduction}
\setcounter{equation}{0}

Moment problems as well as interpolation problems arise in a wide
range of mathematical and physical sciences (see~\cite{A},
\cite{BGM}, \cite{Landau}). The classical power moment problem can
be formulated as follows.

\noindent{\bf The Hamburger moment problem}. Given is a sequence
$\{c_j\}_{j=0}^{\infty}$ of real numbers. Find a positive Borel
measure $d\sigma$ on $\dR$ such that
\[
c_j=\int_{\dR}t^jd\sigma(t),\quad j=0,1,\dots.
\]

In a view of the Hamburger-Nevanlinna theorem (see~\cite{A}), the
moment problem is equivalent to the problem of finding the
Nevanlinna function
$\f(\lambda)\left(=\int_{\dR}\frac{d\sigma(t)}{t-\lambda}\right)$
having the following asymptotic expansions
\[
\f(\lambda)=-\frac{c_{0}}{\lambda}-\frac{c_1}{{\lambda}^2}- \dots
- \frac{c_{2n}}{{\lambda}^{2n+1}}+
o\left(\frac{1}{{\lambda}^{2n+1}}\right) \quad (\lambda=iy,y \to
+\infty)
\]
for all $n\in\dZ_+:=\dN\cup\{0\}$. The latter problem can be
solved by means of the Schur algorithm. This algorithm leads to
the J-fraction
\begin{equation}\label{Jfraction}
-\frac{1}{\lambda-a_0-\displaystyle{\frac{b_0^2}{\lambda-a_1-\displaystyle{\frac{b_1^2}{\ddots}}}}},
\end{equation}
where $a_j$ are real numbers, $b_j$ are positive numbers~\cite{A}
(see also~\cite{GS}). Recall that the theory of J-fractions is
also related to the Pad\'e approximation theory and the theory of
orthogonal polynomials. Under some natural conditions it is
possible to say that all these theories (J-fractions, Pad\'e
approximation and orthogonal polynomials) are equivalent  to one
another. On the other hand, the J-fraction~\eqref{Jfraction}
generates the following infinite Jacobi matrix
\[
\cJ=\left(%
\begin{array}{cccc}
  a_0 & b_0 &  &  \\
  b_0 & a_1 & b_1 &  \\
      & b_1 & a_2 & \ddots \\
      &     & \ddots & \ddots \\
\end{array}%
\right)
\]
In fact, the Jacobi matrix $\cJ$ is a key tool for analyzing the
moment problem as well as the Nevanlinna function $\f$ via
operator methods. For example, using Jacobi matrices techniques
one can prove convergence results for Pad\'e approximants to
Nevanlinna functions (see, for instance,~\cite{BGM}, \cite{Si};
see also~\cite{DD} where generalized Jacobi matrices associated to
indefinite moment problems for generalized Nevanlinna functions
are considered).

The main goal of the present paper is to generalize the
above-described scheme to the case of Nevanlinna-Pick problems in
the class of Nevanlinna functions.  Our approach is based on the
relations between the theory of multipoint Pad\'e approximants
(Pad\'e interpolants), the theory of biorthogonal rational
functions, and generalized eigenvalue problem for two Jacobi
matrices \cite{Zhe_BRF}, \cite{SZ_CMP}, \cite{SZ}, \cite{Zhe_PI}.

%Recall that the ordinary Pad\'e approximation problem is related
%with the theory of orthogonal polynomials. The latter appear also
%as solutions of the ordinary eigenvalue problem for the Jacobi
%(i.e, tri-diagonal) matrix \cite{A}. On the other hand, the
%diagonal Pad\'e approximations give rise to a special type of the
%continued fractions - so-called J-fractions \cite{A}, \cite{BGM},
%\cite{JT}. Under some natural conditions it is possible to say
%that all these theories (J-fractions, Pad\'e approximation and
%orthogonal polynomials) are equivalent one to another.

In theory of biorthogonal rational functions the so-called
continued fractions of $R_{II}$ type appear. These continued
fractions were introduced and studied by Ismail and Masson
\cite{IM}. Nevertheless, note that continued fractions of the same
type were considered earlier in problems connected with rational
interpolation problems (see, e.g. \cite{Wuy}, \cite{DG}). It
appears that the continued fractions of the $R_{II}$ type are
closely connected with the diagonal Pad\'e interpolation problem
from one side and with the theory of generalized eigenvalue
problem for two Jacobi matrices on the other side \cite{Zhe_BRF}.
In turn, both problem are equivalent (under some natural
conditions) to theory of the biorthogonal rational functions (BRF)
\cite{SZ}, \cite{Zhe_PI}. Note that theory of {\it orthogonal}
rational functions studied and developed in \cite{Bulth} can be
considered as a special case of theory of BRF (for details see,
e.g. \cite{Zhe_PI}).

%%%%%%%%%%%%%%%%%%%%%%%%%%%%%%%%%%%%%%%%%%%%%%%%%%%%%%%%%%%%%%%%%%%%%%%%%%%%%%%%%%

\section{Pad\'e interpolation and biorthogonality}
\setcounter{equation}{0} In this section we present basic facts
concerning Pad\'e interpolation and corresponding biorthogonal
rational functions. We follow mostly \cite{Zhe_BRF}, \cite{Zhe_PI}
but some of the result appear to be new.

Let monic polynomials $P_n(z)=z^n + O(z^{n-1})$ satisfy the
$R_{II}$ type recurrence relation \be P_{n+1}(z) + (\alpha_n z +
\beta_n)P_n(z) + r_n (z-a_n)(z-b_n)P_{n-1}(z) =0 \lab{rec_R2} \ee
with initial conditions
$$
P_0=1, \; P_1 = z-\beta_0
$$
Monic property of the polynomials $P_n(z)$ assumes the restriction
upon the recurrence coefficients \be \alpha_0 =-1, \; \alpha_n +
r_n +1 =0, \; n=1,2,\dots \lab{restr_cf} \ee In what follows we
will assume that $r_n \ne 0, \; n=1,2,\dots$ (nondegeneracy).

Introduce the polynomials
$$
A_0=B_0=1, \; A_n(z) = \prod_{k=1}^n(z-a_k), \; B_n(z) =
\prod_{k=1}^n(z-b_k)
$$
As shown by Ismail and Masson \cite{IM} there exists a linear
functional $\sigma$ defined on all rational functions (without a
polynomial part) with the prescribed poles $a_1,b_1,a_2, b_2,
\dots$ by the moments \be c_{nm} =
\sigma\Big\{\frac{1}{A_n(z)B_m(z)}\Big\}, \quad n,m=0,1,2,\dots
\lab{moms_def} \ee such that the orthogonality relation \be
\sigma\Big\{\frac{P_n(z) q_j(z)}{A_n(z)B_n(z)}\Big\} = 0, \quad
j=0,1,\dots, n-1, \lab{ort_IM} \ee holds, where $q_j(z)$ is any
polynomial of degree not exceeding $j$ and
$$
\sigma\Big\{\frac{P_n(z) z^n}{A_n(z)B_n(z)}\Big\} = \kappa_n \ne 0
$$ The normalization coefficients $\kappa_n$ satisfy the
recurrence relation \cite{IM} \be \kappa_{n+1} + \alpha_n \kappa_n
+ r_n \kappa_{n-1}=0 \lab{rec_kappa} \ee It is important to note
that, in contrast to the case of the ordinary orthogonal
polynomials, we can take {\it two} first coefficients $\kappa_0,
\kappa_1$ as arbitrary parameters. Then all further coefficients
$\kappa_2, \kappa_3, \dots$ are determined uniquely through
\re{rec_kappa}.

Note also that if for some $n=n_0>1$ we have
$\kappa_{n_0}=\kappa_{n_0-1}$ then from \re{rec_kappa} and
\re{restr_cf} it follows that $\kappa_{n_0+1}=\kappa_{n_0} =
\kappa_{n_0-1}$ and hence we then have $\kappa_{n}=\kappa_{n_0}$
for all $n \ge n_0-1$. Moreover, we also have from \re{rec_kappa}
$$
r_{n_0-1}(\kappa_{n_0-1} - \kappa_{n_0-2}) =0
$$
Due to our assumption $r_n \ne 0$ we have
$\kappa_{n_0-2}=\kappa_{n_0-1}=\kappa_{n_0}$. Repeating this
process we arrive at condition \be \kappa_1 = \kappa_0.
\lab{deg_kappa_01} \ee We thus have
\begin{proposition}
Condition $\kappa_{n_0} = \kappa_{n_0-1}$ for some $n_0 >1$ is
equivalent to the condition $\kappa_0=\kappa_1$. In this case we
have $\kappa_n \equiv const$ for all $n=0,1,2,\dots$
\end{proposition}
This case will be considered as a degeneration and in what follows
we will assume that $\kappa_1 \ne \kappa_{0}$. Then from this
proposition it follows  $\kappa_n \ne \kappa_{n-1}$ for
$n=2,3,\dots$. Moreover we will assume that $\kappa_n \ne 0$ for
all $n=0,1,2,\dots$.

Introduce the rational functions \cite{Zhe_BRF} \be R_n^{(1)}(z) =
\frac{P_n(z)}{A_n(z)}, \quad R_n^{(2)}(z) = \frac{P_n(z)}{B_n(z)}
\lab{RS_def} \ee It is assumed that zeroes of polynomials $P_n(z)$
do not coincide with points $a_i, b_j$, so rational functions
$R_n^{(1)}(z)$ and $R_n^{(2)}(z)$ have the $[n/n]$ type. Rational
functions $R_n(z)^{(1)}$ have prescribed poles $a_1, a_2, \dots,
a_n$ and rational functions $R_n^{(2)}(z)$ have prescribed poles
$b_1, b_2, \dots, b_n$.

These functions satisfy obvious recurrence relations \be
(z-a_{n+1})R_{n+1}^{(1)}(z) + (\alpha_n z + \beta_n)R_n^{(1)}(z) +
r_n (z-b_n)R_{n-1}^{(1)}(z) =0 \lab{rec_R} \ee

and \be (z-b_{n+1})R_{n+1}^{(2)}(z) + (\alpha_n z +
\beta_n)R_n^{(2)}(z) + r_n (z-a_n)R_{n-1}^{(2)}(z) =0 \lab{rec_S}
\ee On the other hand, these recurrence relations can be rewritten
in terms of the generalized eigenvalue problem (GEVP)
\cite{Zhe_BRF}
$$
J_1 \vec R^{(1)} = z J_2 \vec R^{(1)}
$$
and
$$
J_3 \vec R^{(2)} = z J_2 \vec R^{(2)}
$$
where $\vec R^{(1)}$ is an infinite-dimensional vector with
components $\{ R_0^{(1)}, R_1^{(1)}, \dots \}$ (as well as $\vec
R^{(2)}$) and $J_1, J_2, J_3$ are 3-diagonal (Jacobi) matrices
which entries are obvious from the above recurrence relations for
$R_n^{(1)}, R_n^{(2)}$. As was shown in \cite{Zhe_BRF} the GEVP
leads naturally to theory of biorthogonal rational functions
associated with the polynomials $P_n(z)$ of the $R_{II}$-type.
Here we propose a more simple scheme of construction of the pair
of biorthogonal rational functions.

Introduce the rational functions $U_n(z)$ and $V_n(z)$ by the
formulas: \be U_n(z) = R_n^{(1)}(z) - \xi_n R_{n-1}^{(1)}(z),
\quad V_n(z) = R_n^{(2)}(z) - \xi_n R_{n-1}^{(2)}(z) \lab{UV_def}
\ee where $\xi_n= \kappa_n/\kappa_{n-1}$ (it assumed that
$\xi_0=0$ so that $U_0=V_0=1$). Clearly, the rational functions
$U_n(z)$ have the poles $a_1, a_2, \dots, a_n$ and the rational
functions $V_n(z)$ have the poles $b_1, b_2, \dots, b_n$.

We have
\begin{theorem}\label{orth_T}
The rational functions \re{UV_def} form a biorthogonal system with
respect to the functional $\sigma$: \be \sigma\Big\{ U_n(z) V_m(z)
\Big\} = h_n \: \delta_{nm}, \quad n,m=0,1,\dots \lab{bi_UV} \ee
where the normalization coefficients are
$$
h_n = \frac{\kappa_n}{\kappa_{n-1}}(\kappa_{n-1} - \kappa_n)
$$
\end{theorem}
The proof of this theorem is direct by using orthogonality
relations \re{ort_IM}.

Note that the normalization coefficient is nonzero $h_n \ne 0$ due
to our assumptions on nondegeneracy $\kappa_0 \ne \kappa_1$ and
$\kappa_n \ne 0$.

We can give an equivalent definition of the functions $U_n(z)$ and
$V_n(z)$ using the determinant expressions: \ba &&U_n(z)=
\frac{P_n(a_n)}{\Delta_n} \left |
\begin{array}{cccc} c_{00} & c_{10} & \dots & c_{n,0} \\ c_{01} & c_{11} &
\dots & c_{n,1} \\ \dots & \dots & \dots & \dots\\ c_{0,n-1}&
c_{1,n-1}& \dots & c_{n,n-1}\\ 1& A_1^{-1}(z) & \dots &
A_n^{-1}(z)  \end{array} \right |, \lab{deterU} \ea

\ba &&V_n(z)=\frac{P_n(b_n)}{\Delta_n} \left |
\begin{array}{cccc} c_{00} & c_{01} & \dots & c_{0,n} \\ c_{10} & c_{11} &
\dots & c_{1,n} \\ \dots & \dots & \dots & \dots\\ c_{n-1,0}&
c_{n-1,1}& \dots & c_{n-1,n}\\ 1 & B_1^{-1}(z) & \dots &
B_n^{-1}(z)  \end{array} \right |, \lab{deterV} \ea where \ba
\Delta_n = \left |
\begin{array}{cccc} c_{00} & c_{01} & \dots & c_{0,n-1} \\ c_{10} & c_{11} &
\dots & c_{1,n-1} \\ \dots & \dots & \dots & \dots\\ c_{n-1,0}&
c_{n-1,1} &\dots & c_{n-1,n-1}\end{array} \right |, \lab{Delta}
\ea (It is assumed that $\Delta_0=1$). In what follows we will
assume that
$$
\Delta_n \ne 0, \quad n=1,2,3,\dots
$$
(this is another nondegeneracy condition).

Formulas \re{deterU}, \re{deterV} follow directly from definition
of moments \re{moms_def}. In order to obtain appropriate
coefficients in front of determinantal expressions \re{deterU},
\re{deterV} we can present expression for the rational function
$U_n(z)$ in the following form
$$
U_n(z) = \sum_{k=0}^n \frac{\gamma_{nk}}{A_k(z)}
$$
The leading term in this sum is
$$
\gamma_{nn} = U_n(z) A_n(z)\left|_{z=a_n} \right .
$$
On the other hand we have from the explicit expression \re{UV_def}
$$
U_n(z) A_n(z)\left|_{z=a_n} \right . = P_n(a_n)
$$
whence we obtain the factor $\frac{P_n(a_n)}{\Delta_n}$ in front
of determinant of the formula \re{deterU}. Similarly we obtain the
factor $\frac{P_n(b_n)}{\Delta_n}$ in front of determinant of the
formula \re{deterV}.

Note also that from the determinantal formulas \re{deterU},
\re{deterV} it follows directly that
$$
\sigma\Big\{ U_n(z) V_m(z) \Big\} = \frac{\Delta_{n+1}}{\Delta_n}
P_n(a_n) P_n(b_n) \: \delta_{nm}
$$
Comparing with \re{bi_UV} we obtain an interesting relation \be
P_n(a_n) P_n(b_n) = \frac{\Delta_{n}}{\Delta_{n+1}} \kappa_n (1 -
\kappa_n/\kappa_{n-1}) \lab{PP_kap} \ee From this relation it
follows that condition \be P_n(a_n) P_n(b_n) \ne 0 \lab{PP_ndeg}
\ee guarantees nondegeneracy conditions $\kappa_n \ne 0$,
$\kappa_n \ne \kappa_{n-1}$ and $\Delta_n \ne0$. Thus we will
assume that condition \re{PP_ndeg} holds. It is instructive to
consider what happens when condition \re{PP_ndeg} doesn't hold.
Assume e.g. that $P_n(a_n)=0$ for some $n$. Then the rational
function $R_n^{(1)}(z) = P_n(z)/A_n(z)$ has the order $[n-1/n-1]$,
i.e. it has poles $a_1, a_2, \dots, a_{n-1}$. Corresponding
rational function $U_n(z)$ constructed by \re{UV_def} will also
have poles $a_1, a_2, \dots, a_{n-1}$ which means a degeneration
(absence of the pole $a_n$).

We can present rational functions $U_n(z)$ and $V_n(z)$ in the
form \be U_n(z) = \frac{S_n(z)}{(1-\xi_n)A_n(z)}, \quad V_n(z) =
\frac{T_n(z)}{(1-\xi_n)B_n(z)}, \lab{UV_ST} \ee where $S_n(z)=z^n
+ O(z^{n-1})$ and $T_n(z) = z^n + O(z^{n-1})$ are monic
polynomials of degree $n$. Polynomials $S_n(z), T_n(z)$ are
expressed in terms of polynomials $P_n(z)$: \be S_n(z) =
\frac{P_n(z) - \xi_n (z-a_n) P_{n-1}(z)}{1-\xi_n}, \quad T_n(z) =
\frac{P_n(z) - \xi_n (z-b_n) P_{n-1}(z)}{1-\xi_n} \lab{ST_P} \ee
Moreover $S_0=T_0=1$.

We have
\begin{proposition}
Polynomials $S_n(z), T_n(z)$ satisfy a system of first-order
recurrence relations \ba &&S_{n+1}(z) = \nu_n^{(1)} (z-b_n) S_n(z)
+ \nu_n^{(2)} (z-a_n) T_n(z) \nonumber \\&&T_{n+1}(z) =
\nu_n^{(3)} (z-b_n) S_n(z) + \nu_n^{(4)} (z-a_n) T_n(z), \quad
n=1,2,\dots \lab{rec_ST} \ea where
$$
\nu_n^{(1)}=\frac{\xi_n \beta_n - \xi_n \xi_{n+1} a_{n+1} - r_n
a_n}{r_n(b_n-a_n)}
$$
$$
\nu_n^{(2)}=\frac{\xi_n \beta_n - \xi_n \xi_{n+1} a_{n+1} - r_n
b_n}{r_n(a_n-b_n)}
$$
$$
\nu_n^{(3)}=\frac{\xi_n \beta_n - \xi_n \xi_{n+1} b_{n+1} - r_n
a_n}{r_n(b_n-a_n)}
$$
$$
\nu_n^{(4)}=\frac{\xi_n \beta_n - \xi_n \xi_{n+1} b_{n+1} - r_n
b_n}{r_n(a_n-b_n)}
$$
Note that $\nu_n^{(1)}+ \nu_n^{(2)} = \nu_n^{(3)}+ \nu_n^{(4)}=1$
which is necessary for polynomials $S_{n+1}(z), T_{n+1}(z)$ to be
monic.
\end{proposition}
{\it Proof}. It is sufficient to write down \be U_{n+1}(z) =
\frac{P_{n+1}(z)}{A_{n+1}(z)} - \xi_{n+1} \frac{P_n(z)}{A_n(z)} =
\frac{S_{n+1}(z)}{(1-\xi_{n+1})A_{n+1}} \lab{U_n+1} \ee and
express $P_{n+1}(z)$ in terms of $P_n(z)$ and $P_{n-1}(z)$ using
recurrence relation \re{rec_R2}. This allows one to obtain an
expression of $P_n(z)$ in terms of polynomials $S_n(z),
S_{n+1}(z)$: \be P_n(z) = \zeta_n^{(1)} \left(S_{n+1}(z) -
(z-b_n)S_n(z) \right), \lab{P_S} \ee where
$$
\zeta_n^{(1)} = \frac{r_n(1-\xi_n)}{r_n(b_n-a_{n+1})   -
\xi_n(\beta_n + \alpha_n a_{n+1})}
$$
Analogously \be P_n(z) = \zeta_n^{(2)} \left(T_{n+1}(z) -
(z-a_n)T_n(z) \right), \lab{P_T} \ee where
$$
\zeta_n^{(2)} = \frac{r_n(1-\xi_n)}{r_n(a_n-b_{n+1})   -
\xi_n(\beta_n + \alpha_n b_{n+1})}
$$
Then we arrive at relations \re{rec_ST}.

Note that relations \re{rec_ST} (as well as \re{P_S}, \re{P_T}) do
not hold for $n=0$ because coefficients $\nu_0^{(i)}$ as well as
$a_0, b_0$ are not defined. However, relations \re{rec_ST} will be
valid for $n=0$ if we put
$$
S_0=T_0=1
$$
and \be \xi_0 = \frac{r_0 \kappa_0}{\kappa_0-\kappa_1} \lab{xi0}
\ee whereas $a_0,\: b_0$ and $r_0$ can be {\it arbitrary}
parameters. Indeed it is elementary verified that in this case we
have from relations \re{rec_ST} for $n=0$ \be S_1(z) = z +
\frac{a_1\kappa_1-\beta_0\kappa_0}{\kappa_0-\kappa_1}, \quad
T_1(z) = z + \frac{b_1\kappa_1-\beta_0\kappa_0}{\kappa_0-\kappa_1}
\lab{S1T1} \ee which is compatible with expression for $S_1(z),
T_1(z)$ obtained from \re{ST_P} for $n=1$. Thus we can assume that
relations \re{rec_ST} are valid for $n=0,1,2,\dots$ under
condition \re{xi0}. Note that this condition is formally
equivalent to condition
$$
\kappa_{-1}=\frac{\kappa_0-\kappa_1}{r_0}
$$
obtained from recurrence relation \re{rec_kappa} if one puts $n=0$
(with arbitrary nonzero $r_0$). Equivalently, we can assume that
for $n=0$ coefficients $\nu_n^{(i)}$ take the values
$$
\nu_0^{(1)} = \frac{\kappa_0(\beta_0-a_0)
+\kappa_1(a_0-a_1)}{(b_0-a_0)(\kappa_0-\kappa_1)}, \quad
\nu_0^{(2)} = \frac{\kappa_0(\beta_0-b_0)
+\kappa_1(b_0-a_1)}{(a_0-b_0)(\kappa_0-\kappa_1)}
$$
and
$$
\nu_0^{(3)} = \frac{\kappa_0(\beta_0-a_0)
+\kappa_1(a_0-b_1)}{(b_0-a_0)(\kappa_0-\kappa_1)}, \quad
\nu_0^{(4)} = \frac{\kappa_0(\beta_0-b_0)
+\kappa_1(b_0-b_1)}{(a_0-b_0)(\kappa_0-\kappa_1)}
$$

Vice versa, one can show that starting from the system \re{rec_ST}
with $b_n \ne a_n, \; n=0,1,\dots$ and with initial conditions
$T_0=S_0=1$ one construct a pair of biorthogonal functions
$U_n(z), V_n(z)$ by \re{UV_ST} \cite{SZ}.

The Pad\'e interpolation problem \cite{Becker} (sometimes called
the Cauchy-Jacobi, Newton-Pad\'e or multi-point Pad\'e
approximation problem \cite{Mein}, \cite{BGM}) consists in finding
a pair of polynomials $P_n(z), \: Q_m(z)$ such that \be Y_s
P_n(z_s) = Q_m(z_s), \quad s=0,1,2,\dots, n+m,   \lab{PI_scheme}
\ee where $Y_s$  and $z_s$  are two given complex sequences
($s=0,1,2,\dots$). The degrees of polynomials $P_n(z), \: Q_m(z)$
are less or equal to $n$ and $m$ correspondingly. The rational
function
$$
r_{mn}(z)=\frac{Q_m(z)}{P_n(z)}
$$
is called the rational interpolant corresponding to the sequences
$Y_s$ and $z_s$.

We will consider only the so-called normal case of the Pad\'e
interpolation problem \cite{Becker} meaning that the degrees of
polynomials $P_n(z), Q_m(z)$ are exactly $n$ and $m$ and there are
no common zeros of polynomials $P_n(z)$ and $Q_m(z)$. In the
normal case we have for every pair $(n,m)$ the conditions
\cite{Becker}
$$
r_{m,n+1}(z) \ne r_{mn}(z) \ne r_{m+1,n}(z)
$$
In practice, it is assumed that $Y_s=F(z_s)$ for some given
function $F(z)$ of the complex argument $z$. In this case formula
\re{PI_scheme} gives a rational interpolant $r_{mn}(z)$ of the
function $F(z)$ for the given sequence $z_s$ of interpolation
points. Note that when all interpolation points coincide $z_s=z_0,
s=0,1,2,\dots$, then we have the ordinary Pad\'e approximation in
the point $z_0$. The set $r_{mn}(z), \: m,n=0,1,2,\dots$ is called
the Pad\'e interpolation table for the function $F(z)$.

Consider the so-called diagonal string \cite{SZ}, \cite{Zhe_PI} in
the Pad\'e interpolation table, i.e. the set $r_{n-1,n}(z), \;
n=1,2,\dots$. This means that we are seeking a solution of the
problem \be F(z_s) = \frac{Q_{n-1}(z_s)}{P_n(z_s)}, \quad
s=0,1,2,\dots 2n-1 \lab{diag_PI} \ee Pad\'e interpolants for the
diagonal string satisfy simple orthogonality properties
\cite{Mein}, \cite{SZ} \be [z_0, z_1, \dots z_{2n-1}] \left\{z^j
P_n(z) \right\} = 0, \quad j=0,1,\dots, n-1, \lab{ort_PI} \ee
where \be [z_0, z_1, \dots, z_{2n-1}] \{f(z)\}\equiv \int_{\Gamma}
\frac{f(\zeta) d\zeta }{(\zeta-z_0)(\zeta-z_1)\dots
(\zeta-z_{2n-1})} \lab{Her_DD} \ee is the divided difference of
the order $2n-1$ from the function $f(z)$. It is assumed that the
integration contour $\Gamma$ avoids all singularity points of the
function $f(z)$. Note that formula \re{Her_DD} is called the
Hermite form of the divided difference operation \cite{BGM}.

Orthogonality relation \re{ort_PI} can be extended to
biorthogonality relation for two rational functions $U_n(z),
V_n(z)$ as follows. Consider the diagonal Pad\'e interpolation
problem for the same function $F(z)$ but with slightly modified
interpolation sequence \be F(z_s) = \frac{\tilde Q_{n-1}(z_s)}{\t
P_n(z_s)}, s=0,1,2,\dots, 2n-2, 2n \lab{mod_PI} \ee (i.e. for the
given $n$ we have $2n$ interpolation points as in the previous
scheme \re{diag_PI}, but the final point $z_{2n-1}$ is replaced by
$z_{2n}$). Construct the rational functions \be U_n(z) =
\frac{P_n(z)}{(z-z_1)(z-z_3) \dots (z-z_{2n-1})}, \quad V_n(z) =
\frac{\t P_n(z)}{(z-z_2)(z-z_4) \dots (z-z_{2n})} \lab{UV_PI} \ee
Then the biorthogonality relation \be [z_0, z_1, \dots z_{2n-1}]
\left\{\frac{U_n(z) V_m(z)}{z-z_0} \right\} = h_n \: \delta_{nm},
\quad n,m=0,1,2,\dots \lab{biort_PI} \ee holds with some
normalization constant $h_n \ne 0$ \cite{SZ}, \cite{Zhe_PI}. It is
easily verified that polynomials $P_n(z)$ and $\t P_n(z)$ satisfy
the $R_{II}$ recurrence relations \re{rec_R2} whereas the rational
functions $U_n(z), V_n(z)$ satisfy the generalized eigenvalue
problem of type \re{rec_R}. Thus the generalized eigenvalue
problem for two Jacobi matrices is related with the diagonal
Pad\'e interpolation problem. For further development and
generalizations of this subject see \cite{SZ}, \cite{Zhe_PI},
\cite{Magnus}.

%%%%%%%%%%%%%%%%%%%%%%%%%%%%%%%%%%%%%%%%%%%%%%%%%%%%%%%%%%%%%%%%%%%%%%%%%

\section{Nevanlinna-Pick problems}

In this section we propose a modification of the famous
step-by-step process of solving the Nevanlinna-Pick problem in the
class of Nevanlinna functions~\cite{A},~\cite{AK}.

First, let us recall that a Nevanlinna function is a function
which is holomorphic in the open upper half plane $\dC_+$ and has
a nonnegative imaginary part in $\dC_+$. Let ${\mathbf
N}[\alpha,\beta]$ denote a class of all functions $\f$ having the
representation
\begin{equation}\label{Mar_f}
\f(\lambda)=\int_{\alpha}^{\beta}\frac{d\sigma(t)}{t-\lambda},
\end{equation}
where $d\sigma(t)$ is a finite measure. A function of the class
${\mathbf N}[\alpha,\beta]$ is called a Markov function. Clearly,
a Markov function is also a Nevanlinna function. Moreover, if the
singularities of the Nevanlinna function $\f$ are contained in
$[\alpha,\beta]$ then  $\f\in{\mathbf N}[\alpha,\beta]$ (see, for
instance,~\cite{A}). Let us consider the following Nevanlinna-Pick
problem.

\noindent{\bf Problem NP[$\alpha$,$\beta$].} Given are two
infinite sequences $\{z_k\}_{k=0}^{\infty}$,
$\{w_k\}_{k=0}^{\infty}$ ($z_k\in\dC_+$). Find a function
$\f\in{\mathbf N}[\alpha,\beta]$ such that
\[
\f(z_k)=w_k,\quad k=0,1,2,\dots.
\]

As is known (see~\cite{AK}), the problem $\bf NP[\alpha,\beta]$ is
solvable if and only if the Hermitian forms
\begin{equation}\label{herm_form}
\sum_{j,k=0}^{N}\frac{w_j(z_j-\alpha)-\overline{w}_k(\overline{z}_k-\alpha)
}{z_j-\overline{z}_k}\xi_j\overline{\xi}_k,\quad
\sum_{j,k=0}^{N}\frac{w_j(\beta-z_j)-\overline{w}_k(\beta-\overline{z}_k)
}{z_j-\overline{z}_k}\xi_j\overline{\xi}_k
\end{equation}
are nonnegative definite for all $N\in\dZ_+$.

It is also natural to consider the truncated Nevanlinna-Pick
problem.

\noindent{\bf Problem NP([$\alpha$,$\beta$],n)}. Given are two
finite sequences $\{z_k\}_{k=0}^{n}$, $\{w_k\}_{k=0}^{n}$
($z_k\in\dC_+$). Describe all functions $\f\in{\mathbf
N}[\alpha,\beta]$ satisfying the property
\[
\f(z_k)=w_k,\quad k=0,1,\dots,n.
\]

Note that the problem $\bf NP([\alpha,\beta],n)$ is solvable if
and only if the Hermitian forms~\eqref{herm_form} are nonnegative
definite for $N=0,1,\dots,n$.

The algorithm of solving the Nevanlinna-Pick problems in question
is based on the subsequent statement.

\begin{lemma}\label{schur_step}
Let $\f\in{\mathbf N}[\alpha,\beta]$ and let $z\in\dC_+$ be a
fixed number. Then there exist numbers $a^{(1)},\,a^{(2)}\in\dR$
and $b>0$ such that the function $\tau$ defined by the equality
\begin{equation}\label{flt_step}
\f(\lambda)=-\frac{1}{a^{(2)}\lambda-a^{(1)}+b^2(\lambda-z)(\lambda-\overline{z})\tau(\lambda)}
\end{equation}
belongs to ${\mathbf N}^0[\alpha,\beta]:={\mathbf
N}[\alpha,\beta]\cup\{0\}$.
\end{lemma}

\begin{proof}
Setting ${\displaystyle \Phi(\lambda):=-\frac{1}{\f(\lambda)}}$,
define the function
\begin{equation}\label{leq1}
\Psi(\lambda)=\frac{\Phi(\lambda)-\Phi(z)}{\Phi(\lambda)-\overline{\Phi(z)}}:\frac{\lambda-z}{\lambda-\overline{z}}.
\end{equation}
Due to the Schwartz lemma, we have that
\[
|\Psi(\lambda)|\le 1, \quad \Im\lambda>0.
\]
So, the function $\widetilde{\Psi}_1$ defined from the relation
\begin{equation}\label{leq2}
\Psi(\lambda)=\frac{\widetilde{\Psi}_1(\lambda)-i}{\widetilde{\Psi}_1(\lambda)+i}
\end{equation}
is a Nevanlinna function. Plugging~\eqref{leq2} into~\eqref{leq1},
one obtains
\begin{equation}\label{leq3}
\widetilde{\Psi}_1(\lambda)=-i\frac{
\Phi(\lambda)(2\lambda-z-\overline{z})-\overline{\Phi(z)}(\lambda-z)-\Phi(z)(\lambda-\overline{z})
} {
\Phi(\lambda)(z-\overline{z})-\Phi(z)(\lambda-\overline{z})+\overline{\Phi(z)}(\lambda-z)
}.
\end{equation}
Now, let us consider the following function
\[
{\Psi}_1(\lambda):=\widetilde{\Psi}_1(\lambda)+i\frac{2\lambda-(z+\overline{z})}{z-\overline{z}}=
\widetilde{\Psi}_1(\lambda)+\frac{2\lambda-(z+\overline{z})}{2\Im
z}.
\]
Obviously, $\Psi_1$ is a Nevanlinna function. Taking into
account~\eqref{leq3}, $\Psi_1$ admits the following representation
\begin{equation}\label{leq4}
\Psi_1(\lambda)=-\frac{\Im \Phi(z)}{\Im
z}\frac{(\lambda-z)(\lambda-\overline{z})}{\Phi(\lambda)-\frac{\Im\Phi(z)}{\Im
z}\lambda+\frac{\Im \Phi(z)\overline{z}}{\Im z}}.
\end{equation}
Finally, introducing
\[
\tau(\lambda)=-\frac{1}{\Psi_1(\lambda)}\in{\mathbf N},
b=\frac{\Im \Phi(z)}{\Im z}>0, a^{(2)}=\frac{\Im\Phi(z)}{\Im
z}\in\dR, a^{(1)}=-\frac{\Im \Phi(z)\overline{z}}{\Im z}\in\dR,
\]
one can easily transform~\eqref{leq4} into~\eqref{flt_step}. To
complete the proof, it is sufficient to observe that, due
to~\eqref{Mar_f} and~\eqref{flt_step}, all singularities of $\tau$
are contained in $[\alpha,\beta]$.
\end{proof}

\begin{remark} The transformation~\eqref{flt_step} could be
viewed as a substitute for the Schwartz lemma. A similar
to~\eqref{flt_step} transformation for Caratheodory functions was
proposed in~\cite{DG}.
\end{remark}

\begin{remark}
Substituting $\lambda$ for $z$ and $\overline{z}$
in~\eqref{flt_step} we get
\[
\f(z)=-\frac{1}{a^{(2)}z-a^{(1)}},\quad
\f(\overline{z})=-\frac{1}{a^{(2)}\overline{z}-a^{(1)}}.
\]
Expressing from the above relations $a^{(1)}$ and $a^{(2)}$, one
can obtain the following formulas
\begin{equation}\label{first_step_f}
a^{(2)}=-\frac{\Im\frac{1}{\f(z)}}{\Im z},\quad
a^{(1)}=-\frac{\Im\frac{1}{\f(z)}}{\Im z}z+\frac{1}{\f(z)}.
\end{equation}
It is easy to see that the numbers $a^{(1)}$, $a^{(2)}$ are
uniquely determined by~\eqref{first_step_f}. Further,
equality~\eqref{flt_step} can be rewritten as follows
\begin{equation}\label{flt_t}
b^2\tau(\lambda)=-\frac{\frac{1}{\f(\lambda)}+a^{(2)}\lambda-a^{(1)}}{(\lambda-z)(\lambda-\overline{z})}.
\end{equation}
In fact, the number $b$ can be chosen arbitrary. So, to be
definite we always choose $b>0$ in the following way
\[
b^2=\int_{\alpha}^{\beta}d\sigma(t).
\]
In this case, the function $\tau$ possesses the integral
representation~\eqref{Mar_f} with a probability measure.
\end{remark}
\begin{remark}
It also easily follows from the theory of generalized Nevanlinna
functions (see~\cite{DLLS}, \cite{DHS01}, \cite{KL77}) that the
right-hand side of~\eqref{flt_t} is a Nevanlinna function.
\end{remark}
\begin{remark}
By comparing the first terms in asymptotic expansions of the right
hand side and left hand side of~\eqref{flt_step}, we see that
\begin{equation}\label{connection}
a^{(2)}=\left(\int_{\alpha}^{\beta}d\sigma(t)\right)^{-1}+b^2.
\end{equation}
\end{remark}

Now, we are in a position to solve the problem $\bf
NP([\alpha,\beta],n)$. Let the given problem $\bf
NP([\alpha,\beta],n)$ be solvable and let $\f$ be a solution of
the problem $\bf NP([\alpha,\beta],n)$. Due to
Lemma~\ref{schur_step}, $\f_0:=\f$ admits the following
representation
\begin{equation}\label{gSchur}
\f(\lambda)=-\frac{1}{a^{(2)}_0\lambda-a^{(1)}_0+b^2_0(\lambda-z_0)(\lambda-\overline{z}_0)\f_1(\lambda)},
\end{equation}
where $\f_1\in{\mathbf N}^0[\alpha,\beta]$. From~\eqref{gSchur} we
see that
\[
\f_1(\lambda)=-\frac{\frac{1}{\f(\lambda)}+a_{0}^{(2)}\lambda-a_0^{(1)}}{b_0^2(\lambda-z_0)(\lambda-\overline{z}_0)}.
\]
So, if $\f_1\not\equiv 0$ then it is a solution of the problem
$\bf NP([\alpha,\beta],n-1)$ with the sequences
$\{z_k\}_{k=1}^{n}$ and $\{w_k^{(1)}\}_{k=1}^{n}$, where
\[
w_k^{(1)}=\f_1(z_k)=-\frac{\frac{1}{w_k}+a_{0}^{(2)}z_k-a_0^{(1)}}{b_0^2(z_k-z_0)(z_k-\overline{z}_0)}.
\]
Therefore, the original problem $\bf NP([\alpha,\beta],n)$ is
reduced to the problem $\bf NP([\alpha,\beta],n-1)$. Similarly,
the problem $\bf NP([\alpha,\beta],n-1)$ can be reduced to the
problem $\bf NP([\alpha,\beta],n-2)$ and so on. Finally, one has a
sequence of the linear fractional transformations
\[
\f_{j}(\lambda)=-\frac{1}{a^{(2)}_j\lambda-a^{(1)}_j+b^2_j(\lambda-z_j)(\lambda-\overline{z}_j)\f_{j+1}(\lambda)}\quad(j=0,1,\dots,n)
\]
having the following matrix representations
\begin{equation}\label{W_j}
\cW_j(\lambda)=\begin{pmatrix} 0               & \displaystyle{-\frac{1}{b_j(\lambda-\overline{z}_j)}}\\
                               b_j(\lambda-z_j)& \displaystyle{\frac{a^{(2)}_j\lambda-a^{(1)}_j}{b_j(\lambda-\overline{z}_j)}}
                            \end{pmatrix}
                            \quad (j=0,1,\dots,n).
\end{equation}
If the above-described algorithm consists of exactly $n+1$ steps
then we say that the problem $\bf NP([\alpha,\beta],n)$ is {\it
nondegenerate}. So, we have proved the following theorem which
gives the complete solution of the problem $\bf
NP([\alpha,\beta],n)$.

\begin{theorem}[\cite{A}] Any solution $\f$ of the nondegenerate
problem $\bf NP([\alpha,\beta],n)$ admits the following
representation
\begin{equation}\label{NPTrunc}
    \f(\lambda)=\frac{w_{11}(\lambda)\tau(\lambda)+w_{12}(\lambda)}
    {w_{21}(\lambda)\tau(\lambda)+w_{22}(\lambda)},
\end{equation}
where $\tau\in{\mathbf N}^0[\alpha,\beta]$ and
\begin{equation}\label{W}
\cW_{[0,n]}(\lambda)=(w_{ij}(\lambda))_{i,j=1}^2:=\cW_0(\lambda)\cW_1(\lambda)\dots
\cW_n(\lambda).
\end{equation}
\end{theorem}

It should be also remarked that $\cW_j$ is the Blaschke-Potapov
factor~\cite{ADRS}, \cite{Po}.

%%%%%%%%%%%%%%%%%%%%%%%%%%%%%%%%%%%%%%%%%%%%%%%%%%%%%%%%%%%%%%%%%%%%%%%%%%%%%%%%%%%%%%%%%%%%%%%%%%%%%%%%%%%%%%%%%%%

\section{Rational systems related to Nevanlinna-Pick problems}

Let us suppose that the given Markov function  has the integral
representation~\eqref{Mar_f} with a probability measure $d\sigma$
which support contains infinitely many points, i.e.
\[
\int_{\alpha}^{\beta}d\sigma(t)=1.
\]
 Let us also suppose that for the given sequence
$\{z_k\}_{k=0}^{\infty}$ there exists $\delta>0$ such that
\begin{equation}\label{zcond}
\Im z_k>\delta,\quad k=0,1,2,\dots.
\end{equation}
In this case, it follows from the uniqueness theorem for analytic
functions that the numbers $z_k$ and $w_k:=\f(z_k)$ ($k\in\dZ_+$)
uniquely determine the function $\f$. So, the Nevanlinna-Pick
problem with the data $\{z_k\}_{k=0}^{\infty}$,
$\{w_k\}_{k=0}^{\infty}$ has a unique solution.

Since $\f$ is not rational the given data give rise to infinitely
many steps of the step-by-step process. So, we have infinitely
many linear fractional transformations of the
form~\eqref{flt_step} which lead to the following continued
fraction
%\begin{equation}
\begin{multline}\label{ContF}
-\frac{1}{a^{(2)}_0\lambda-a_0^{(1)}-\displaystyle{\frac{b_0^2(\lambda-z_0)(\lambda-\overline{z}_0)}{a^{(2)}_1\lambda-a_1^{(1)}-
\displaystyle{\frac{b_1^2(\lambda-z_1)(\lambda-\overline{z}_1)}{\ddots}}}}}\\
=-\frac{1}{a^{(2)}_0\lambda-a_0^{(1)}}
\begin{array}{l} \\ -\end{array}
\frac{b_0^2(\lambda-z_0)(\lambda-\overline{z}_0)}{a^{(2)}_1\lambda-a_1^{(1)}}
\begin{array}{ccc} \\ - \end{array}
\frac{b_{1}^2(\lambda-z_1)(\lambda-\overline{z}_1)}{a^{(2)}_2\lambda-a_2^{(1)}}\begin{array}{l} \\
-\cdots \end{array}.
\end{multline}
%\end{equation}
The continued fraction~\eqref{ContF} is an $R_{II}$-fraction
(see~\cite{IM}). Consider the $(n+1)$-th convergent of the
continued fraction~\eqref{ContF}
\[
R_n(\lambda):=-\frac{1}{a^{(2)}_0\lambda-a_0^{(1)}}
\begin{array}{l} \\ -\end{array}
\frac{b_0^2(\lambda-z_0)(\lambda-\overline{z}_0)}{a^{(2)}_1\lambda-a_1^{(1)}}
\begin{array}{ccc} \\ - \end{array}
\begin{array}{l} \\
\cdots \end{array}\begin{array}{ccc} \\ - \end{array}
\frac{b_{n-1}^2(\lambda-z_{n-1})(\lambda-\overline{z}_{n-1})}{a^{(2)}_n\lambda-a_n^{(1)}}.
\]
It is obvious that $R_n$ is a solution of the problem $\bf
NP([\alpha,\beta],n)$, i.e. the following equality holds true
\[
R_n(z_k)=w_k=\f(z_k),\quad k=0,1,\dots,n.
\]
\begin{definition}\label{PadeInt}
The $[L/M]$ multipoint Pad\'e approximant  for a function $\f$ at
the points $\{\alpha_k\}_{k=1}^{\infty}$ is defined as a ratio
\[
f^{[L/M]}(\lambda)=\frac{A^{[L/M]}(\lambda)}{B^{[L/M]}(\lambda)}
\]
of two polynomials $A^{[L/M]}$, $B^{[L/M]}$ of formal degree $L$
and $M$, respectively, such that
\[
f^{[L/M]}(\alpha_k)=\f(\alpha_k),\quad k=1,\dots,L+M+1.
\]
\end{definition}

Since $R_n$ is real, the rational function $R_n$ is the $[n/n]$
multipoint Pad\'e approximant for $\f$ at the points
$\{\infty,z_0,\overline{z}_0,\dots,z_n,\overline{z}_n,\dots\}$.

It is well-known that to every continued fraction there
corresponds a recurrence relation. In particular, for the
continued fraction~\eqref{ContF} a recurrence relation takes the
following form
\begin{equation}\label{rec_rel}
u_{j+1}-(a_j^{(2)}\lambda-a_j^{(1)})u_{j}+b_{j-1}^2(\lambda-z_{j-1})(\lambda-\overline{z}_{j-1})u_{j-1}=0\quad(j\in\dN).
\end{equation}
Define polynomials of the first kind $P_j(\lambda)$ as solutions
$u_j=P_j(\lambda)$ of the system~\eqref{rec_rel} with the initial
conditions
\begin{equation}\label{InConP}
u_{0}=1,\quad u_1=a_0^{(2)}\lambda-a_0^{(1)}.
\end{equation}
Similarly, the polynomials of the second kind $Q_j(\lambda)$ are
defined as solutions $u_j=Q_j(\lambda)$ of the
system~\eqref{rec_rel} subject to the following initial conditions
\begin{equation}\label{InConP2}
u_{0}=0,\quad u_1=-1.
\end{equation}

Note that in our setting~\eqref{ort_IM} is transformed into the
following orthogonality relations (see also~\cite{GL})
\begin{equation}\label{NP_ort}
\int_{\alpha}^{\beta}t^jP_{n+1}(t)\frac{d\sigma(t)}{|t-z_0|^2\dots|t-z_{n}|^2}=0,\quad
j=0,1,\dots, n.
\end{equation}

 It follows from the theory of continued fractions that
$R_n(\lambda)={\displaystyle\frac{Q_{n+1}(\lambda)}{P_{n+1}(\lambda)}}$
(see, for details,~\cite{JT}).

The recurrence relation~\eqref{rec_rel} can be renormalized to the
following one
\begin{equation}\label{r_rl}
b_{j}(z_{j}-\lambda)\widehat{u}_{j+1}-(a_j^{(2)}\lambda-a_j^{(1)})\widehat{u}_{j}+
b_{j-1}(\overline{z}_{j-1}-\lambda)\widehat{u}_{j-1}=0\quad(j\in\dN),
\end{equation}
where
\[
\widehat{u}_0=u_0,\quad
\widehat{u}_j={\displaystyle\frac{u_j}{b_0\dots
b_{j-1}(z_0-\lambda)\dots(z_{j-1}-\lambda)}}\quad (j\in\dN).
\]
Relation~\eqref{NP_ort} implies that
\begin{equation}\label{NP_ort_2}
\int_{\alpha}^{\beta}\widehat{P}_{n+1}(t)\frac{1}{t-\overline{z}_j}d\sigma(t)=0,\quad
j=0,\dots,n.
\end{equation}
Now, setting
\[
\xi_0=0,\quad
\xi_j=\left(\int_{\alpha}^{\beta}\frac{t^{j+1}d\sigma(t)}{|t-z_0|^2\dots|t-z_{j}|^2}\right)
\left(\int_{\alpha}^{\beta}\frac{t^{j}d\sigma(t)}{|t-z_0|^2\dots|t-z_{j-1}|^2}\right)^{-1}\quad(j\in\dN)
\]
one can see that the simple linear combinations
$\widehat{P}_{j}-\xi_j\widehat{P}_{j-1}$ ($j\in\dZ_+$) give
orthogonalization of the system
\[
\left\{1, \frac{1}{\lambda-z_0}, \frac{1}{\lambda-z_1},\dots
\right\}
\]
of rational functions (see~Theorem~\ref{orth_T}, see
also~\cite{BEZ}). It should be also noted here that systems of
orthogonal rational functions related to Nevanlinna-Pick problems
were proposed in~\cite{BL},~\cite{LL},~\cite{Nja} (see
also~\cite{Bulth}).

 Further, relation~\eqref{r_rl} can be rewritten as
follows
\begin{equation}\label{gen_ee}
z_jb_{j}\widehat{u}_{j+1}+a_j^{(1)}\widehat{u}_j+z_{j-1}b_{j-1}\widehat{u}_{j-1}=
\lambda(b_{j}\widehat{u}_{j+1}+a_j^{(2)}\widehat{u}_j+b_{j-1}\widehat{u}_{j-1})\quad
(j\in\dN).
\end{equation}
The system~\eqref{gen_ee} gives us the possibility to rewrite the
Cauchy problem~\eqref{rec_rel},~\eqref{InConP} in the matrix form
\[
J_{[0,\infty)}^{(1)}\pi(\lambda)=\lambda
J_{[0,\infty)}^{(2)}\pi(\lambda),
\]
where $\pi(\lambda)=\left(\widehat{P}_0(\lambda),
\widehat{P}_1(\lambda),\dots,\widehat{P}_j(\lambda),\dots
\right)^{\top}$ and
\[
J_{[0,\infty)}^{(1)}=\left(%
\begin{array}{cccc}
  a_0^{(1)} & \overline{z}_0b_0 &  &  \\
  z_0b_0 & a_1^{(1)} & \overline{z}_1b_1 &  \\
      & z_1b_1 & a_2^{(1)} & \ddots \\
      &     & \ddots & \ddots \\
\end{array}%
\right)\quad
J_{[0,\infty)}^{(2)}=\left(%
\begin{array}{cccc}
  a_0^{(2)} & b_0 &  &  \\
  b_0 & a_1^{(2)} & b_1 &  \\
      & b_1 & a_2^{(2)} & \ddots \\
      &     & \ddots & \ddots \\
\end{array}%
\right).
\]
We denote by $\ell^2_{[0,n]}$ the space of $(n+1)$ vectors with
the usual inner product.  Define a standard basis in
$\ell^2_{[0,n]}$ by the equalities
\[
e_{j}=\{{\delta_{l,k}\}}_{k=0}^{n},\quad j=0,1,\dots,n.
\]
Let $J_{[j,k]}^{(1)}$ ($J_{[j,k]}^{(2)}$) be a submatrix of
$J^{(1)}_{[0,\infty)}$ ($J^{(2)}_{[0,\infty)}$), corresponding to
the linear subspace spanned by the vectors $e_{l},\dots,e_k$
$(0\le j\le k \le n)$, that is,
\[
J^{(1)}_{[j,k]}=
\begin{pmatrix}
  a_j^{(1)} & \overline{z}_jb_j & {\bf 0} \\
  z_jb_j & \ddots &    \\
      {\bf 0}&  & a_k^{(1)} \\
      \end{pmatrix},\quad
      J^{(2)}_{[j,k]}=
\begin{pmatrix}
  a_j^{(2)} & b_j &  {\bf 0}  \\
  b_j & \ddots &   \\
  {\bf 0}&  & a_k^{(2)} \\
    \end{pmatrix}.
\]

\begin{proposition}
The matrix $J_{[0,n]}^{(2)}$ is positive definite for all
$n\in\dZ_+$.
\end{proposition}
\begin{proof}
Let us consider the Hermitian form
\begin{equation}\label{jform2}
\left(J_{[0,n]}^{(2)}\xi,\xi\right)=a_0^{(2)}|\xi_0|^2+b_0\overline{\xi}_0\xi_1+
b_0\xi_0\overline{\xi}_1+a_1^{(2)}|\xi_1|^2+\dots+a_n^{(2)}|\xi_n|^2.
\end{equation}
Due to~\eqref{connection} and our assumptions, we have that
$a_j^{(2)}=1+b_j^2$.
%Now, the statement of proposition immediately
%follows from the following factorization
%\[
%J_{[0,n]}^{(2)}=
%\begin{pmatrix}
%a_{0}^{(2)}&b_{0}&&{\bf 0}\\
%b_{0}&\ddots&\ddots&\\
%&\ddots&a_{n-1}^{(2)}&b_{n-1}\\
%{\bf 0}&&b_{n-1}&a_{n}^{(2)}\\
%\end{pmatrix}=
%\begin{pmatrix}
%1 &b_{0}&&{\bf 0}\\
%0&\ddots&\ddots&\\
%&\ddots&1&b_{n-1}\\
%{\bf 0}&&0&1\\
%\end{pmatrix}
%\begin{pmatrix}
%1&{0}&&{\bf 0}\\
%b_{0}&\ddots&\ddots&\\
%&\ddots&1&0\\
%{\bf 0}&&b_{n-1}&1\\
%\end{pmatrix}.
%\]
Therefore, one can
rewrite the form~\eqref{jform2} in the following manner
\[
\left(J_{[0,n]}^{(2)}\xi,\xi\right)=|\xi_0|^2+|b_0\xi_0+\xi_1|^2+\dots+|b_{n-1}\xi_{n-1}+\xi_{n}|^2+|b_n\xi_n|^2\ge
0.
\]
Thus, $J_{[0,n]}^{(2)}$ is a positive definite matrix.
\end{proof}

Finally, we should note that for the matrix $J_{[0,\infty)}^{(2)}$
the following factorization holds true
\[
J_{[0,\infty)}^{(2)}=\left(%
\begin{array}{cccc}
  a_0^{(2)} & b_0 &  &  \\
  b_0 & a_1^{(2)} & b_1 &  \\
      & b_1 & a_2^{(2)} & \ddots \\
      &     & \ddots & \ddots \\
\end{array}%
\right)=
\left(%
\begin{array}{cccc}
  1 & b_0 &  &  \\
  0 & 1 & b_1 &  \\
      & 0 & 1 & \ddots \\
      &     & \ddots & \ddots \\
\end{array}%
\right)
\left(%
\begin{array}{cccc}
  1 & 0 &  &  \\
  b_0 & 1 &0 &  \\
      & b_1 & 1 & \ddots \\
      &     & \ddots & \ddots \\
\end{array}%
\right).
\]

\section{$m$-functions of linear pencils}

In this section we give a matrix representation of multipoint
Pad\'e approximants for Markov functions.
\begin{definition}
The function
\begin{equation}\label{Weyl1}
m_{[j,n]}(\lambda)=\left((J_{[j,n]}^{(1)}-\lambda
J_{[j,n]}^{(2)})^{-1}e_j,e_j\right)%\quad e_0=(1,0,\dots,0)^{\top}
\end{equation}
will be called the $m$-function of the linear pencil
$J_{[j,n]}^{(1)}-\lambda J_{[j,n]}^{(2)}$.
\end{definition}

To see the correctness of the above definition it is sufficient to
rewrite~\eqref{Weyl1} in the following form
\begin{equation}\label{cor_Weyl}
m_{[j,n]}(\lambda)=\left((J_{[j,n]}^{(2)})^{-1}(J_{[j,n]}^{(1)}(J_{[j,n]}^{(2)})^{-1}-\lambda)^{-1}e_j,e_j\right).
\end{equation}
From~\eqref{cor_Weyl} one can conclude that $m_{[j,n]}$ is a
Nevanlinna function.

\begin{proposition}
The $m$-functions $m_{[j,n]}$ and $m_{[j+1,n]}$ are related by the
equality
\begin{equation}\label{Riccati}
m_{[j,n]}=-\frac{1}{a^{(2)}_j\lambda-a^{(1)}_j+b^2_j(\lambda-z_j)(\lambda-\overline{z}_j)m_{[j+1,n]}(\lambda)}.
\end{equation}
\end{proposition}
\begin{proof}
Consider the following block representation of the matrix
$J_{[j,n]}^{(1)}-\lambda J_{[j,n]}^{(2)}$
\[
J_{[j,n]}^{(1)}-\lambda J_{[j,n]}^{(2)}=\begin{pmatrix}
a^{(1)}_j-a^{(2)}_j\lambda &  B\\
B^{*}      &   J_{[j+1,n]}^{(1)}-\lambda J_{[j+1,n]}^{(2)}\\
\end{pmatrix},
\]
where $B=(b_j(\overline{z}_j-\lambda),0,\dots, 0)$. According to
the Frobenius formula~\cite[Section~0.7.3]{HJ} the matrix
$(J_{[j,n]}^{(1)}-\lambda J_{[j,n]}^{(2)})^{-1}$ has the following
block representation
\begin{equation}\label{Frobenius}
(J_{[j,n]}^{(1)}-\lambda J_{[j,n]}^{(2)})^{-1}=\begin{pmatrix}
\left(a^{(1)}_j-a^{(2)}_j\lambda-B^*(J_{[j+1,n]}^{(1)}-\lambda J_{[j+1,n]}^{(2)})^{-1}B\right)^{-1}&*\,\,\\
* &*\,\,\\
\end{pmatrix}.
\end{equation}
Plugging~\eqref{Frobenius} into \eqref{Weyl1}, one
obtains~\eqref{Riccati}.
\end{proof}
\begin{corollary}\label{zeros}
The following equalities hold true
\begin{equation}\label{mformulas}
m_{[0,n]}(\lambda)=R_n(\lambda)=\frac{Q_{n+1}(\lambda)}{P_{n+1}(\lambda)}\in\mathbf{
N}[\alpha,\beta].
\end{equation}
\end{corollary}
\begin{proof}
Relation~\eqref{Riccati} implies that the rational functions
$m_{[0,n]}$ and $R_n$ have the same expansions into
$R_{II}$-fractions. So, $m_{[0,n]}=R_n$. By using standard
argumentation, from~\eqref{NP_ort} one can conclude that {\it all
the zeros of $P_{n+1}$ are contained in $[\alpha,\beta]$}
(see~\cite{A},~\cite{GL}). The latter means that the Nevanlinna
function $m_{[0,n]}$ belongs to $\mathbf{N}[\alpha,\beta]$.
\end{proof}
So, now one can say that $R_n$ is a solution of $\bf{
NP([\alpha,\beta],n)}$. By using standard argumentation,
from~\eqref{mformulas} we can conclude the following result.
\begin{corollary}\label{interlace}
The zeros of $P_{n+1}$ and $Q_{n+1}$ are interlace.
\end{corollary}

Below, we will need the following statement.
\begin{corollary}\label{spectr}
The spectrum
$\sigma\left(J_{[0,n]}^{(1)}(J_{[0,n]}^{(2)})^{-1}\right)$ of the
matrix $J_{[0,n]}^{(1)}(J_{[0,n]}^{(2)})^{-1}$ is contained in
$[\alpha,\beta]$.
\end{corollary}
\begin{proof} From the formula for calculation of inverse
matrices,~\eqref{Weyl1}, and~\eqref{mformulas} one can see that
\[
m_{[0,n]}(\lambda)=\frac{\det(J_{[1,n]}^{(1)}-\lambda
J_{[1,n]}^{(2)})}{\det(J_{[0,n]}^{(1)}(J_{[0,n]}^{(2)})^{-1}-\lambda)\det(J_{[0,n]}^{(2)})}=
\frac{Q_{n+1}(\lambda)}{P_{n+1}(\lambda)}.
\]
So, the statement immediately follows from
Corollary~\ref{interlace} and the fact that all the zeros of
$P_{n+1}$ are contained in $[\alpha,\beta]$.
\end{proof}

\begin{remark} It should be remarked that, for the case of the Laurent orthogonal
polynomials,
 a similar scheme with two matrices and $m$-functions were considered in~\cite{CKA}.
\end{remark}
%%%%%%%%%%%%%%%%%%%%%%%%%%%%%%%%%%%%%%%%%%%%%%%%%%%%%%%%%%%%%%%%%%%%%%%%%%%%%%%%%%%%%%%%%%%%%%%%%%%%%%%%

\section{A convergence result for multipoint Pad\'e approximants}

The goal of this section is to prove an analog of Markov's
convergence theorem by making use of the operator representation
of multipoint Pad\'e approximants.

We begin with an auxiliary statement.
\begin{lemma}\label{lem_in}
The following inequalities hold true
\[
\left((J_{[0,n]}^{(2)})^{-1}e_0,e_0\right)\le 1\quad (n\in\dZ_+).
\]
\end{lemma}
\begin{proof}
The proof is by induction. First, note that
\[
\left((J_{[n,n]}^{(2)})^{-1}e_n,e_n\right)=\frac{1}{a_n^{(2)}}=\frac{1}{1+b_n^2}\le
1\quad (n\in\dZ_+).
\]
Suppose that
$\displaystyle{\left((J_{[k+1,n]}^{(2)})^{-1}e_{k+1},e_{k+1}\right)}\le
1$. It follows from the Riccati equation~\cite[formula~(2.15)]{GS}
(see also~\eqref{Riccati}) that
\[\begin{split}
\left((J_{[k,n]}^{(2)})^{-1}e_k,e_k\right)=&\frac{1}{a_k^{(2)}-b_k^2\left((J_{[k+1,n]}^{(2)})^{-1}e_{k+1},e_{k+1}\right)}\\
=&\frac{1}{1+b_k^{2}-b_k^2\left((J_{[k+1,n]}^{(2)})^{-1}e_{k+1},e_{k+1}\right)}\le
1.\end{split}
\]
\end{proof}

Now, we are ready to prove the main result of this section.
\begin{theorem}[\cite{GL}]\label{MarkovMP}
Let $\f\in{\mathbf N}[\alpha,\beta]$ and let the sequence
$\{z_k\}_{k=1}^{\infty}$ satisfy the condition~\eqref{zcond}. Then
the sequence $f^{[n/n]}=R_n$ converges to $\f$ locally uniformly
in $\dC\setminus[\alpha,\beta]$.
\end{theorem}
\begin{proof}
We first recall the well-known estimate for the resolvent of
self-adjoint operator $J$ (for instance
see~\cite[Theorem~V.3.2]{Kato})
\begin{equation}\label{estimat}
\Vert(J-\lambda)^{-1}\Vert\le\frac{1}{\dist(\lambda,\sigma(J))}.
\end{equation}
Next, observe that the operator
$J_{[0,n]}^{(1)}(J_{[0,n]}^{(2)})^{-1}$ is self-adjoint with
respect to the following inner product
\[
\left((J_{[0,n]}^{(2)})^{-1}x,y\right)\quad x,y\in\dC^{n+1}.
\]
 Taking into account the
representation~\eqref{mformulas},~\eqref{cor_Weyl}, the
Cauchy-Schwartz inequality,~\eqref{estimat},
Corollary~\ref{spectr}, and Lemma~\ref{lem_in}, we obtain
\begin{equation}\label{main_est}
\begin{split}
|R_{n}(\lambda)|&=\left|\left(J_{[0,n]}^{(2)})^{-1}(J_{[0,n]}^{(1)}(J_{[0,n]}^{(2)})^{-1}-\lambda)^{-1}e_0,e_0\right)\right|\\
&\le\frac{\left((J_{[0,n]}^{(2)})^{-1}e_0,e_0\right)}{\dist(\lambda,[\alpha,\beta])}\le
\frac{1}{\dist(\lambda,[\alpha,\beta])}.
\end{split}
\end{equation}
It follows from~\eqref{main_est} and Montel's theorem that the
family $\{R_n\}$ is precompact in the topology of locally uniform
convergence in $\dC\setminus[\alpha,\beta]$. Note that
\[
R_n(z_k)=\f(z_k),\quad n\ge k.
\]
Thus, applying the Vitali theorem completes the proof.
\end{proof}

\begin{remark} Theorem~\ref{MarkovMP} was proved in~\cite{GL} by
means of another method. The rates of convergence of multipoint
Pad\'e approximants was also given in~\cite{GL}. The operator
interpretation of the rates of convergence and a more detailed
analysis of the underlying linear pencil will be given in the
forthcoming paper.
\end{remark}

\noindent{\bf Acknowledgments}. MD expresses his gratitude to
Professor V.A.~Derkach for carefully reading the manuscript and
giving many helpful comments. The authors also thank the referees
for helpful suggestions.

%%%%%%%%%%%%%%%%%%%%%%%%%%%%%%%%%%%%%%%%%%%%%%%%%%%%%%%%%%%%%%%%%%%%%%%%%%%%%%%%%%%%%%%%%%%%%%%%%%%%%%%

\end{document}